\def\Bbb#1{{\bf #1}}

\def\fnote#1{\footnote}
\def\blacksquare{\hbox{\vrule width 4pt height 4pt depth 0pt}}

\def\cwdash{\relbar\joinrel}

\def\cwleftpar#1#2{\leftskip #1 \rightskip #2 plus 1fill}
\def\cwrightpar#1#2{\leftskip #1 plus 1fill \rightskip #2}
\def\cwcenterpar#1#2{\leftskip #1 plus 1fill \rightskip #2 plus 1fill}
\def\cwfullpar#1#2{\leftskip#1\rightskip#2}

\def\cwoutdent#1#2{\llap{\hbox to #1{#2 \hss}}\ignorespaces}
\def\cwparbegin#1#2#3#4#5{
	\ifcase #1 \cwleftpar{#2}{#3}
	\or \cwrightpar{#2}{#3}
	\or \cwcenterpar{#2}{#3}
	\else \cwfullpar{#2}{#3}\fi
	\ifcase #4 \baselineskip = 1.5\baselineskip
	\or \baselineskip = 2\baselineskip
	\or \baselineskip = 3\baselineskip
	\else \baselineskip = 1\baselineskip\fi
	\ifdim #5 > 0in \else \noindent \fi
	\noindent\ignorespaces}
\documentclass{article}
\begin{document}
\advance \vsize by -1\baselineskip
\def\makefootline{
{\vskip \baselineskip \noindent \folio                                  \par
}}

\vspace{2ex}
\noindent {\Huge Special Bases for Derivations of\\[0.4ex]
			 Tensor Algebras}\\[0.3ex]
\noindent {\Large III. Case along Smooth Maps\\[0.5ex]
		with Separable Points of Selfintersection}\\

\vspace{2ex}

\noindent Bozhidar Z. Iliev
\fnote{0}{\noindent $^{\hbox{}}$Permanent address:
Laboratory of Mathematical Modeling in Physics,
Institute for Nuclear Research and \mbox{Nuclear} Energy,
Bulgarian Academy of Sciences,
Boul.\ Tzarigradsko chauss\'ee~72, 1784 Sofia, Bulgaria\\
\indent E-mail address: bozho@inrne.bas.bg\\
\indent URL: http://theo.inrne.bas.bg/$^\sim$bozho/}

\vspace{3ex}

{\bf \noindent Published: Communication JINR, E5-92-543, Dubna, 1992}\\[1ex]
\hphantom{\bf Published: }
http://www.arXiv.org e-Print archive No.~math.DG/0305061\\[2ex]

2003 PACS numbers: 02.40.Ma, 02.40.Vh,  04.20.Cv, 04.90.+h\\[2ex]

\noindent
{\small
The \LaTeXe\ source file of this paper was produced by converting a
ChiWriter 3.16 source file into
ChiWriter 4.0 file and then converting the latter file into a
\LaTeX\ 2.09 source file, which was manually edited for correcting numerous
errors and for improving the appearance of the text.  As a result of this
procedure, some errors in the text may exist.
}\\[2ex]

	\begin{abstract}
Necessary and/or sufficient conditions are studied for the existence,
uniqueness and holonomicity of bases in which on sufficiently general subsets
of a differentiable manifold the components of derivations of the tensor
algebra over it vanish. The linear connections and the equivalence principle
are considered form that point of view.
	\end{abstract}\vspace{3ex}

 {\bf I. INTRODUCTION}

\medskip
 In connection with the equivalence principle [1], as well as from purely
 mathematical reasons $[2-5]$, of importance is the problem for existence of
 local (holonomic or anholonomic [2]) coordinates (bases) in which the
 components of a linear connection [3] vanish on  some  subset, usually
 submanifolds, of a differentiable manifold [3]. This problem has been solved
 for torsion free, i.e. symmetric, linear connections $[3, 4]$ in the cases
 at a point $[2-5]$, along a smooth path without selfintersections $[2, 5]$
 and on a neighborhood $[2, 5]$. These results were generalized in our
 previous works $[6, 7]$ for arbitrary, with or without torsion, derivations
 of the tensor algebra over a given differentiable manifold [3] and, in
 particular, for arbitrary linear connections. General results of this kind
 can be found in [8], where a criteria for the existence of the
 above-mentioned special bases (coordinates) on submanifolds of a space with
 a symmetric affine connection is given.

 The present work generalizes the results from $[6-8]$ and deals with the
 problems for existence, uniqueness and holonomicity of special bases in
 which the components of a derivation of the tensor algebra over a
 differentiable manifold vanish on some its subset of a sufficiently general
 type (Sect. II and III). In particular, this derivation may be a linear
 connection (Sect. IV). In this context we also make conclusions concerning
 the general validity and the mathematical formulation of the equivalence
 principle in a class of gravity theories (Sect. $V)$.

 For further reference purposes, as well as for the exact statement of the
 above problems, we below reproduce below a few simple facts about
 derivations of a tensor algebras that can be found in $[6, 7]$ or derived
 from the ones in [3].

 Every  derivation  $D$  of  the tensor algebra over a differentiable
 manifold  $M$ admits a unique representation in the form (see  $[3], ch.I$,
 proposition $3.3) D=L_{X}+S$, in which $L_{X}$ is the Lie derivative along
 the vector field $X$ and $S$ is tensor field of type (1,1) considered here
 as a differentiation [3]. Both $X$ and $S$ are uniquely defined by D.  If
 $S$  is  a  map  from  the set of $C^{1}$ vector fields into the  tensor
 fields of type (1,1) and $S:X\mapsto S_{X}$, then the equation
 \[
 D^{S}_{X}=L_{X}+S_{X}\qquad (1)
 \]
   defines a derivation of the tensor algebra
 over $M$ for any $C^{1}$ vector field $X [3]$. As the map $S$ will hereafter
 be assumed fixed, such a derivation, i.e. $D^{S}_{X}$, will be called an
 $S$-derivation along $X$ and will be denoted for brevity by $D_{X}$. An
 $S$-derivation is a map $D$ such that $D:X\mapsto D_{X}$ where $D_{X}$ is an
 $S$-derivation along X.

 Let $\{E_{i}, i=1,\ldots  ,n:=\dim(M)\}$ be a (coordinate or noncoordinate
 [2]) local basis of vector fields in the tangent to $M$ bundle. It is
 (an)holonomic if the vectors $E_{1}, \ldots  , E_{n}$ (do not) commute [2].
  The  local  components  $(W_{X})^i_{.j}$ of $D$  with respect to $\{E_{i}\}$
 are defined by the equation
\[
  D_{X}(E_{j})=(W_{X})^{i}_{.j}E_{i}. \qquad(2)
 \]
 and their explicit form is
 \[
(W_{X})^{i}_{.j}=(S_{X})^{i}_{.j}-E_{j}(X^{i})+C^{i}_{.jk}X^{k},\qquad (3)
 \]
 where  $X(f)$  denotes  the  action of $X=X^{k}E_{k}$ on  $C^{1}$functions
 $f$, i.e.  $X(f):=X^{k}E_{k}(f)$,  and  $C^{i}_{.jk}$define the commutators
 of the basic vectors, i.e.
 \[
  [E_{j},E_{k}]=C^{i}_{.jk}E_{i}\qquad (4)
 \]
  If we make the change
$\{E_{i}\} \to \{E_{i^\prime }=A^{i}_{i^\prime }E_{i}\}$, where
$A:= [A^{i}_{i^\prime }]=:[ A^{i^\prime }_{i}] ^{-1}$is a nondegenerate
matrix function, then from (2) we can see that $(W_{X})^{i}_{.j}$transform
into
\[
  (W_{X})^{i^\prime }_{..j^\prime }=A^{i^\prime }_{i}A^{j}_{j^\prime
 }(W_{X})^{i}+A^{i^\prime }_{i}X(A^{i}_{j^\prime }),\qquad (5)
 \]
 which, if we  introduce  the  matrices  $W :=  (W_{X})^{i}_{.j}  $  and
 $W^\prime :=  (W_{X})^{i^\prime }_{..j^\prime }  $, will read
 \[
  W=A^{-1}[W_{X}A+X(A)],\qquad (5^\prime )
 \]
  where the superscript  is understood   a  first  matrix  index and
 $X(A):=X^{k}E_{k}(A)=  X^{k}E_{k}(A^{i}_{i^\prime })  $.

If $\nabla $ is a linear connection with local components $\Gamma
 ^{i}_{.jk}($see e.g. $[2-4])$, then is fulfilled $[2, 3]$
 \[
  \nabla _{X}(E_{j})=(\Gamma ^{i}_{.jk}X^{k})E_{i}.\qquad (6)
 \]
  Hence,  comparing (2) and (6), we see that $D_{X}$is a covariant
 differentiation along $X$ iff
 \[
  (W_{X})^{i}_{.j}=\Gamma ^{i}_{.jk}X^{k}\qquad (7)
 \]
 for some functions $\Gamma ^{i}_{.jk}$.

 Let  $D$  be  an $S$-derivation and $X$ and $Y$ be vector fields. The
 torsion operator $T$  of $D$ is defined as
\[
  T^{D}(X,Y):=D_{X}Y-D_{Y}X-[X,Y].\qquad (8)
 \]
 The $S$-derivation $D$ is called torsion free if $T^{D}=0$.

 For a linear connection $\nabla $, due to (7), we have
 \[
 (T^{\nabla }(X,Y))^{i}=T^{i}_{.jk}X^{j}Y^{k},\qquad (9)
 \]
  where $[2, 3]$
 $T^{i}_{.kl}:=-(\Gamma ^{i}_{.kl}-\Gamma ^{i}_{.lk})-C^{i}_{.kl}$ are the
 components of the torsion tensor of $\nabla $.

 Further  we  shall  investigate  the problem for existence of special  bases
 $\{E_{i^\prime }\}$ in which $W$  =0 for an $S$-derivation $D$ along any or
 fixed vector field X. Hence, due to $(5^\prime )$, we shall have to solve
 the  equation  $W_{X}(A)+X(A)=0$ with respect to A under conditions that will be presented below.

 \newpage

 {\bf II. DERIVATIONS ALONG ARBITRARY VECTOR FIELDS}

\medskip
 This  section is devoted to the existence and some properties  of  special
 bases $\{E_{i^\prime }\}$, defined in a neighborhood of a subset $U$ of  the
 manifold $M$, in which the components of an $S$-derivation $D$ along  an
 {\it arbitrary} vector field $X$ vanish on U.

 {\bf Proposition 1.} If for some $S$-derivation $D$ there exists a basis
 $\{E_{i^\prime }\}$ in which $W$  $\mid _{U}=0$ for every vector field $X$,
 then $D$ is a linear connection on the set $U\subset $M.

  {\bf Remark.}  On  the set $U\subset M$ the derivation $D$ is a linear
 connection  if  $(cf. (7))$  in  some,  and  hence  in  any, basis
 $\{E_{i}\}$ is fulfilled
 \[
  W_{X}(x)=\Gamma _{k}(x)X^{k}(x),\qquad (10)
 \]
  where  $x\in U, X=X^{k}E_{k}$and $\Gamma _{k}$ are some matrix functions
 on $U$. (Evidently, a  linear connection on $M$ is also on $U a$ linear
 connection for every $U$; see (7)).

  {\it Proof.}  If  we  fix  a  basis $\{E_{i}\}$ and $E_{i^\prime
 }=A^{i}_{i^\prime }E_{i}$, then by the  definition of $\{E_{i^\prime }\}$,
 we have $W$  $\mid _{U}=0$, i.e. $W$  $(x)=0$ for $x\in U$, which in
 conformity  with  $(5^\prime )$  is  equivalent to (10) with $\Gamma
 _{k}=-(E_{k}(A))A^{-1}$, $A= [ A^{i}_{i^\prime }]  .\blacksquare $

 The opposite statement to proposition 1 is generally not true  and for its
 exact formulation we shall need some preliminary results and explanations.

 Let  $p$  be an integer, $p\ge 1$, and the Greek indices $\alpha $ and
 $\beta $ run from 1 to p. Let $J^{p}$be a neighborhood in ${\Bbb R}^{p}$ and
 $\{s^{\alpha }\}=\{s^{1},\ldots  ,s^{p}\}$ be  (Cartesian) coordinates in
 ${\Bbb R}^{p}$.

  {\bf Lemma 2.} Let $Z_{\alpha }:J^{p}\cwdash \to GL(m,{\Bbb R}), GL(m,{\Bbb
 R})$ being the group of $m\times m$ non-degenerate matrices  on  ${\Bbb R}$,
 be $a C^{1}$  matrix-valued function on $J^{p}$. Then the  initial-value
 problem
 \[
 \frac{\partial Y}{\partial s^a}\Big|_s
 =Z_{\alpha }(s)Y, \quad
 Y|_{s=s_{0}}={\Bbb I}:=  [\delta  ^{i}_{j}]^{m}_{i,j=1},\qquad (11)
\]
 in  which ${\Bbb I}$  is  the  unit  matrix of the corresponding size, $s\in
J^{p}$, $s_{0}\in J^{p}$  is fixed and $Y$ is $m\times m$ matrix function on
$J^{p}$, has a solution,  denoted by $Y=Y(s,s_{0};Z_{1},\ldots  ,Z_{p})$,
which is unique and smoothly depends on all its arguments, if and only if
\[
  R_{\alpha \beta }(Z_{1},\ldots  ,Z_{p}):=\partial Z_{\alpha }/\partial
s^{\beta }- \partial Z_{\beta }/\partial s^{\alpha }+ Z_{\alpha }Z_{\beta }-
Z_{\beta }Z_{\alpha }=0.\qquad (12)
\]

{\it Proof.}  According  to the results  in $[9], ch$. VI the integrability
conditions for (11) are
\[
  0 = \partial ^{2}Y/\partial s^{\alpha }\partial s^{\beta }- \partial
 ^{2}Y/\partial s^{\beta }\partial s^{\alpha }= \partial (Z_{\beta
 }Y)\backslash \partial s^{\alpha }- \partial (Z_{\alpha }Y)\backslash
 \partial s^{\beta }=
 \]
 \[
 = (\partial Z_{\beta }\backslash \partial s^{\alpha })Y - (\partial
 Z_{\alpha }\backslash \partial s^{\beta })Y + Z_{\beta }Z_{\alpha }Y -
 Z_{\alpha }Z_{\beta }Y = -R_{\alpha \beta }(Z_{1},\ldots  ,Z_{p})Y.
 \]
  Hence (11) has a unique solution iff (12) is satisfied.\blacksquare

 Let  $p\le n:=\dim(M), \alpha ,\beta =1,\ldots  ,p$ and $\mu ,\nu
=p+1,\ldots  $,n. Let $\gamma :J \to M$  be a $C^{1}$  map. We shall suppose
 that for any $s\in J^{p}$ there exists its  $(p$-dimensional)  neighborhood
 $J_{s}\subset J^{p}, s\in J$  such that the restricted  map  $\gamma \mid
 _{J_{s}}:J_{s}\cwdash \to M$ is without selfintersections, i.e. in
 $J_{s}$ does not  exist  points  $s_{1}$  and  $s_{2}\neq s_{1}$  with  the
 property $\gamma (s_{1})=\gamma (s_{2})$. This  assumption   is   equivalent
 to   the  one  that  the  points  of  selfintersections  of $\gamma $, if
  any, can be separated by neighborhoods. With  $J^{p}_{s}$we shall denote
 the union on all neighborhoods $J_{s}$ with the  above  property;  evidently
 $J^{p}_{s}$is the maximal neighborhood of $s$ in  which $\gamma $ is without
 selfintersections.

 Let  at  first  suppose  $J^{p}_{s}=J^{p}$,  i.e. $\gamma $ to be without
 selfintersection,  and  that  $\gamma (J^{p})$  be  contained in only one
 coordinate  neighborhood $V$ of M.

 Let  us  fix  some one-to-one onto $C^{1}$ map $\eta :J^{p}\times
 J^{n-p}\mapsto M$ such that  $\eta (\cdot ,\mathbf{ t}_{0})=\gamma $  for  a
 fixed  $\mathbf{ t}_{0}\in J^{n-p}$,  i.e.  $\eta (s,\mathbf{ t}_{0})=\gamma (s),
 s\in J^{p}$.  In  $V\cap \eta (J^{p},J^{n-p})$  we  define  coordinates
 $\{x^{i}\}$ by putting $(x^{1}(\eta (s,\mathbf{ t}))$,  $\ldots  ,x^{n}(\eta
 (s,\mathbf{ t}))):=(s,\mathbf{ t})\in {\Bbb R}^{n}, s\in J^{p}, \mathbf{ t}\in J^{n-p}$.

 {\bf Proposition 3.} Let $\gamma :J^{p}  \to M$ be $C^{1}$, without
 selfintersections and  $\gamma (J^{p})$  lies  in only one coordinate
 neighborhood. Let on $\gamma (J^{p})$ the  derivation  $D$  be a linear
 connection. Then there exists a defined  in  a  neighborhood of $\gamma
 (J^{p})$ basis $\{E_{i^\prime }\}$ in which the components  of $D$ along
 every vector field vanish on $\gamma (J^{p})$ if and only if in the
 above-defined coordinates $\{x^{i}\}$ is fulfilled
 \[
  [R_{\alpha \beta }(-\Gamma _{1}\circ \gamma ,\ldots  ,-\Gamma _{p}\circ
 \gamma )]|_{J^{p}}=0, \alpha ,\beta =1,\ldots  ,p,\qquad (13)
  \]
  where
 $R_{\alpha \beta }(\ldots  )$  is  defined  by (12) for $m=n$ and
 $(s^{1},\ldots  ,s^{p})=s\in J^{p}$,
 i.e.
 \[
   [R_{\alpha \beta }(-\Gamma _{1}\circ \gamma ,\ldots  ,-\Gamma _{p}\circ
 \gamma )](s)=\partial \Gamma _{\alpha }(\gamma (s))/\partial s^{\beta }-
 \partial \Gamma _{\beta }(\gamma (s))/\partial s^{\alpha }+
 \]
 \[
 + (\Gamma _{\alpha }\Gamma _{\beta }- \Gamma _{\beta }\Gamma _{\alpha
 })|_{\gamma (s)}.\qquad (14)
 \]

  {\bf Remark.} In the case when $D$ is a
 symmetric affine connection this result was obtained by means of another
 method in [8].

{\it Proof.} The following considerations will be done in the above-defined
neighborhood $V\cap \eta (J^{p},J^{n-p})$ and coordinates $\{x^{i}\}$. Let
 $E_{i}=\frac{\partial}{\partial x^i}$.
We shall look for a basis $\{E_{i^\prime }=A^{i}_{i^\prime
}E_{i}\}$ in which $W$  $(\gamma (s))=0, s\in J^{p}$. By $eq. (5^\prime )$
the existence of $\{E_{i^\prime }\}$ is equivalent to  the existence of the
$A=  A^{i}_{i^\prime }  $, transforming $\{E_{i}\}$ into $\{E_{i^\prime }\}$,
and such that $[A^{-1}(W_{X}A+X(A))]\mid _{\gamma (s)}=0$ for every X. But as
on $\gamma (J^{p}) D$ is a linear connection, the $eq. (10)$ is valid for
some matrix-valued functions $\Gamma _{k}$and $x\in \gamma (J^{p})$,
consequently $A$ must be a solution of $\Gamma $  $_{^\prime }(x)=0$, i.e. of
\[
(\Gamma _{k}(\gamma (s)) A(\gamma (s))
 + (\partial A/\partial x^{k})| _{\gamma (s)}=0,\quad s\in J^{p}.\qquad (15)
\]

 By expanding $A(\eta (,\mathbf{ t})), s\in J^{p}, \mathbf{ t}\in J^{n-p}$into a
power series  with respect to $(\mathbf{ t}-\mathbf{ t}_{0})$ it can be shown that
(15) has a solution if  and only if the integrability conditions (13)  are
valid.  Besides,  if (13) take place, than the general solution of (15) is
\[
  A(\eta (s,\mathbf{ t}))
 =\Bigl\{
 {\Bbb I}-\sum_{\lambda =p+1}^n\Gamma _{\lambda }(\gamma (s))
  [x^{\lambda }(\eta (s,\mathbf{ t}))-x^{\lambda }(\gamma (s))]
  \Bigr\} \times
  \]
  \[
  \times    Y(s,s_{0};-\Gamma _{1}\circ \gamma ,\ldots  ,-\Gamma _{p}\circ
  \gamma )B_{0} + \sum_{\mu ,\nu=p+1}^n B_{\mu \nu }(s,\mathbf{ t};\eta
  )[x^{\mu }(\eta (s,\mathbf{ t}))-
  \]
  \[
   -x^{\mu }(\gamma (s))][x^{\nu }(\eta (s,\mathbf{ t}))-x^{\nu }(\gamma
  (s))],\qquad (16)
 \]
   where $s_{0}\in J^{p}$ and the nondegenerate matrix
  $B_{0}$ are fixed and the matrices $B_{\mu ,\nu }, \mu ,\nu =p+1,\ldots
  ,n$, together with their  derivatives, are bounded when $\mathbf{ t}\cwdash \to
  \mathbf{ t}_{0}$. (The fact that into (16) enter only sums from $p+1$ to $n$ is
  a consequence from $x^{\alpha }(\eta (s,\mathbf{ t}))=x^{\alpha }(\gamma
  (s))=s^{\alpha }$, i.e. $x^{\alpha }(\eta (s,\mathbf{ t}))-x^{\alpha }(\eta
  (s,\mathbf{ t}_{0}))=x^{\alpha }(\eta (s,\mathbf{ t}))-x^{\alpha }(\gamma
  (s))\equiv 0, \alpha =1,\ldots  $,p.)

  Thus  bases  $\{E_{i^\prime }\}$ in which $W$  =0 exist iff (13) is
  satisfied.  If (13) is valid, then the bases $\{E_{i^\prime }\}$ are
  obtained from $\{E_{i}=\partial/\partial x^{i}\}$ by means of linear
 transformations the matrices of which must have the form (16).\blacksquare

 Now we are ready to consider a general smooth $(C^{1})$ map $\gamma
  :J^{p}\cwdash \to M$ whose points of selfintersection, if any, can be
  separated by neighborhoods. Let for any $r\in J^{p}$be chosen a coordinate
  neighborhood $V_{\gamma (r)}$of $\gamma (r)$ in M. Let there be fixed $a
  C^{1}$one-to-one onto map $\eta _{r}:J^{p}_{r}\times J^{n-p}\to M$
 such that $\gamma|_{J^{p}_{r}}=\eta _{r}(\cdot ,\mathbf{ t}^{r}_{0})$ for
 some $\mathbf{ t}^{r}_{0}\in J^{n-p}$. In the neighborhood $V_{\gamma (r)}\cap
  \eta _{r}(J^{p}_{r},J^{n-p})$ of $\gamma (J^{p}_{r})\cap V_{\gamma (r)}$we
  introduce local coordinates $\{x^{i}_{r}\}$ defined by $(x^{1}_{r}(\eta
  _{r}(s,\mathbf{ t})),\ldots  ,x^{n}_{r}(\eta _{r}(s,\mathbf{ t}))):= :=(s,\mathbf{
  t})\in {\Bbb R}^{n}$, where $s\in J^{p}_{r}$and $\mathbf{ t}\in J^{n-p}$are
  such that $\eta _{r}(s,\mathbf{ t})\in V_{\gamma (r)}$.

 {\bf Theorem 4.} Let the points of selfintersection of the $C^{1}$map
  $\gamma :J^{p}\cwdash \to M$, if any, be separable by neighborhoods and let
  on $\gamma (J^{p})$ the $S$-derivation $D$ be a linear connection, i.e.
  $eq. (10)$ to be valid. Then in some neighborhood of $\gamma (J^{p})$
  exists a basis $\{E_{i^\prime }\}$ in which the components of $D$ along
  every vector field vanish on $\gamma (J^{p})$ if and only if for every
  $r\in J$ in the above defined local coordinates $\{x^{i}_{r}\}$ is
  fulfilled
 \[
  [R_{\alpha \beta }(-\Gamma _{1}\circ \gamma ,\ldots  ,-\Gamma _{p}\circ
  \gamma )](s)=0, \alpha ,\beta =1,\ldots  ,p,\qquad (17)
  \]
 where  $\Gamma_{\alpha }$ are calculated by means of (10) in $\{x^{i}_{r}\},
  R_{\alpha \beta }$ are given by  (14) and $s\in J^{p}_{r}$ is such that
  $\gamma (s)\in V_{\gamma (r)}$.

 {\it Proof.}  For any $r\in J^{p}$, the restricted map $\gamma|_{^\prime
  J^{p}_{r}}:^\prime J^{p}_{r}\to M$, where  $^\prime
  J^{p}_{r}:=\{s\in J^{p}_{r}, \gamma (s)\in V_{\gamma (s)}\}$,  is  without
  selfintersections  (see the  above  definition  of $J^{p}_{r})$
 and $\gamma|_{^\prime J^{p}_{r}}(^\prime J^{p}_{r})=\gamma (^\prime
  J^{p}_{r})$ lies in the coordinate neighborhood $V_{\gamma (r)}$.

  	So,  if    a  basis $\{E_{i^\prime }\}$  with the described property
  exists,  then, by proposition 3, eqs. (17) are identically satisfied.

  	On  the   opposite, if  (17)  are  valid,  then,  again  by
  proposition  3, for  every $r\in J^{p}$in a neighborhood $^\prime V_{r}$ of
  $\gamma (^\prime J^{p}_{r})$ in $V_{\gamma (r)}$  exists  a  basis
 $\{E^{r}_{i^\prime }\}$  in which the components of $D_{X}$ along  every
  vector  field $X$ vanish on $\gamma (^\prime J^{p}_{r})$. From the
  neighborhoods $^\prime V_{r}$  we  can construct  a  neighborhood  $V$  of
  $\gamma (J^{p})$, e.g. we can  put
 $V=\cup_{r\in J^p} ^\prime V_{r}$, but,
  generally, $V$ is sufficient to be taken as  a  union  of $^\prime V_{r}$ or
  some, but not all $r\in J^{p}$.  On $V$ we can obtain a basis
  $\{E_{i^\prime }\}$  with  the  needed property  by putting $E_{i^\prime
  }\mid _{x}=E^{r}_{i^\prime }\mid _{x}$if $x$ belongs to  only  one
  neighborhood $^\prime V_{r}$and if $x$ belongs to more than one
  neighborhood  $^\prime V_{r}$we can choose $\{E_{i^\prime }\mid _{x}\}$ to
  be the basis $\{E^{r}_{i^\prime }\mid x\}$ for some arbitrary  fixed  $r$
  with this property. (Note that generally the so obtained  basis  is not
  continuous in the regions containing intersections of several neighborhoods
  $^\prime V_{r}.)\blacksquare $

 {\bf Proposition 5.} If on the set $U\subset M$ there exists bases in which
  the components of some $S$-derivation along every vector field vanish on
  $U$, then all of them are obtained from one another by linear
  transformations whose coefficients are such that the action on them of the
  corresponding basic vectors vanishes on U.

 {\it Proof.}The proposition is a simple corollary from $(5^\prime
  ).\blacksquare $

{\bf Proposition 6.} If for some $S$-derivation $D$ there exists a locally
  holonomic basis in which the components of $D$ along every vector field
  vanish on the set $U\subset M$, then $D$ is torsion free on U. On the
  opposite, if $D$ is torsion free on $U$ and bases in which the components
  of $D$ along every vector field vanish on $U$ exist, then all of them are
  holonomic on $U$, i.e. their basic vectors commute on U.

 {\it Proof.}  If $\{E_{i^\prime }\}$ is a basis with the mentioned
  property, i.e. $W (x)=0$  for  every $X$ and $x\in U$, then using (2) and
  (8) (see also eq.  (15)  from [6]), we fined $T^{D}(E_{i^\prime
  },E_{j^\prime })\mid _{U}=-[E_{i^\prime },E_{j^\prime }]\mid _{U}$and
  consequently  $\{E_{i^\prime }\}$  is holonomic on $U$, i.e. $[E_{i^\prime
  },E_{j^\prime }]\mid _{U}=0$, iff $0=T^{D}(X,Y)\mid _{U}=\{X^{i^\prime
  }Y^{j^\prime }T^{D}(E_{i^\prime },E_{j^\prime })\}\mid _{U}$for  every
  vector fields $X$ and $Y$, which is  equivalent to $T^{D}\mid _{U}=0$.

 On  the opposite, let $T^{D}\mid _{U}=0$. We want to prove that any basis
  $\{E_{i^\prime }\}$  in which $W$  =0 is holonomic on U. The holonomicity
  on $U$ means  $0=[E_{i^\prime },E_{j^\prime }]\mid _{U}=\{-A^{k^\prime
  }_{k}(E_{j^\prime }(A^{k}_{i^\prime })-E_{i^\prime }(A^{k}_{j^\prime
  }))E_{k^\prime }\}\mid _{U}$.  But (see proposition  1  and  (10)) the
  existence of $\{E_{i^\prime }\}$ is equivalent to$_{U}W_{X}\mid
  _{U}=(\Gamma _{k}X^{k})\mid $  for  some  functions $\Gamma _{k}$and every
  X. These two facts, combined with  (2)  and  (8),  show  that  $(\Gamma
  _{k})^{i}_{.j}=(\Gamma _{j})^{i}_{.k}$.  Using  this  and  $\{\Gamma
  _{k}A+ \partial A/\partial x^{k}\}\mid _{U}=0 ($see the proof of
  proposition 1), we find  $E_{j^\prime }(A^{k}_{i^\prime })\mid
  _{U}=-A^{j}_{j^\prime }A^{i}_{i^\prime }(\Gamma _{j})^{k}_{.i}\mid
  _{U}=E_{i^\prime }(A^{k}_{j^\prime })\mid _{U}$and therefore $[E_{i^\prime
  },E_{j^\prime }]\mid _{U}=0 ($see above), i.e. $\{E_{i^\prime }\}$ is
  holonomic on U.\blacksquare

\medskip
\medskip
 {\bf III. DERIVATIONS ALONG FIXED VECTOR FIELD}

\medskip
	 As  it  was said in our previous works $[6, 7]$ the problem for
  existence and the properties of special bases for derivations along  {\it
  fixed}  vector  field  is not very interesting from the viewpoint of  its
  applications. By that reason we shall briefly outline only some  results
  concerning it. The following two propositions are almost evident (cf. resp.
  propositions 1 and 5).

  {\bf Proposition  7.}  If  for  the  $S$-derivation  $D_{X}$along a fixed
  vector  field  $X$ exists a basis $\{E_{i^\prime }\}$ in which the
  components of $D_{X}$  vanish  on the set $U\subset M$, then on $U$
 $D_{X}$ is a covariant differentiation  along $X$, i.e. for the given $X$
  the $eq.  (10)$ is valid on U.

{\bf Proposition 8.} If on the set $U\subset M$ there exists bases in which
  the components of an $S$-derivation along a fixed vector field vanish, then
  all of them are obtained from one another by linear transformations, the
  matrices of which are such that the action of $X$ on them vanishes on U.

	 The existence of special bases in which the components of $D_{X}$,
  with  a  {\it fixed} $X$, vanish on some set $U\subset M$ significantly
  differs from  the  same  problem  for  $D_{X}$  with an {\it arbitrary} $X
  ($see Sect. II). In  fact,  if  $\{E_{i^\prime }=A^{i}_{i^\prime }E_{i}\},
  \{E_{i}\}$  being  a fixed basis on $U$, is such a  basis  on $U$, i.e. $W$
  $\mid _{U}=0$, then, due to $(5^\prime )$, its existence is equivalent  to
  the  one  of  $A:=  A^{i}_{i^\prime }  $ for which $(W_{X}A+X(A))\mid
  _{U}=0$ for the {\it given}  X.  As  $X$ is fixed, the values of A at two
  different points,  say $x,y\in U$, are connected through the last equation
  if and only if $x$  and  $y$ lie on one and the same integral path (curve)
  of $X$, the part  of which between $x$ and $y$ belongs entirely to U.
  Hence, if $\gamma :J\cwdash \to M$,  $J$  being  an ${\Bbb R}$-interval, is
  (a part of) an integral path of $X$, i.e.  at  $\gamma (s), s\in J$ the
  tangent to $\gamma $ vector field    is
  $\dot\gamma(s):=X| _{\gamma (s)}$, then
  along  $\gamma $  the equation $(W_{X}A+X(A))\mid _{U}=0$ reduces to
  $dA/ds|_{\gamma (s)}=\dot\gamma(A)|_{s}
  =(X(A))\mid _{\gamma (s)}=-W_{X}(\gamma (s))A(\gamma (s))$.
  The general solution of this equation is
  \[
   A(s;\gamma )=Y(s,s_{0};-W_{X}\circ \gamma )B(\gamma ),\qquad (18)
  \]
   where
  $s_{0}\in J$ is fixed, $Y=Y(s,s_{0};Z)$, with $Z$ being $a C^{1}$ matrix
  function  of $s$, is the unique solution of the initial-value problem (see
  [9], ch. IV, \S1)
  \[
  dY/ds=ZY, Y|_{s=s_{0}}={\Bbb I},\qquad (19)
  \]
  and  the nondegenerate matrix $B(\gamma )$ may depend only on $\gamma $,
  but not on  s. (Note that (19) is a special case of (11) for $p=1$ and by
  lemma 2  it has always a unique solution because of $R_{11}(Z_{1})\equiv 0$
  for $p=1.)$

 From the above considerations follows

 {\bf Proposition 9.} For any $S$-derivation along a fixed vector field on
  every set $U\subset M$ there exist bases, in which the components of that
  derivation vanish on U.

\medskip
\medskip
 {\bf IV. LINEAR CONNECTIONS}

\medskip
 The  results  of Sect. II can directly be applied to the case  of  linear
  connections.  As this is more or less trivial, we shall  present without
  proofs only three such consequences.

 {\bf Corollary  10.}
 Let  the points of selfintersection of the $C^{1}$ map
$\gamma :J^{p} \to  M$, if any, be separable by neighborhoods, $\nabla $ be a
linear connection on  $M$  with local components $\Gamma ^{i}_{.jk}$ (in a
basis $\{E_{i}\})$ and $\Gamma _{k}:= [\Gamma ^{i}_{.jk}]^{n}_{i,j=1}$.  Then
in a neighborhood of $\gamma (J^{p})$ exists a basis  $\{E_{i^\prime }\}$ in
which the components of $\nabla $ vanish on $\gamma (J^{p})$, i.e. $\Gamma
_{k^\prime }\mid _{\gamma (J^{p_{)}}}=0$,  iff  for  every  $r\in J^{p}$  in
the coordinates $\{x^{i}_{r}\}$ (defined before theorem 4) is satisfied (17)
in which $\Gamma _{\alpha }, \alpha =1,\ldots  ,p$ are part of the components
of $\nabla $ in $\{x^{i}_{r}\}$ and $s\in J^{p}$ is such that $\gamma (s)\in
V_{\gamma (r)}$.

 {\bf Corollary 11.} If on the set $U\subset M$ there exist bases in which
  the components of a linear connection vanish on $U$, then all of these
  bases are obtained from one another by linear transformations, the matrices
  of which are such that the action of the corresponding basic vectors on
  them vanishes on U.

{\bf Corollary 12.} Let for some linear connection on a neighborhood of a
  set $U\subset M$ exist local continues bases in which the connection's
  components vanish on U. Then one, and hence any, such basis is holonomic on
  $U$ iff the connection is torsion free on U.

\medskip
\medskip
 {\bf V. CONCLUDING REMARKS}

\medskip
	 As the main result of this work is expressed by theorem 4, we shall
  make some comments on it. First of all, it expresses a sufficiently general
  necessary and sufficient condition for existence of the considered here
special bases for derivations, in particular linear connections. For
instance, it covers that problem on arbitrary submanifolds. In this sense,
its special cases are the results in our previous papers $[6, 7]$.

 If  $p=0$  or  $p=1$,  then  the  condition  (17)  is identically
satisfied,  i.e. $R_{\alpha \beta }=0 ($see (14)). Hence in these two cases
special  bases,  we are searching for, always exist (respectively at a point
 or along a path), which was already established in [6] and [7] respectively.

	 In  the  other  limiting case, $p=n:=\dim(M)$, it is easily seen
 that  the  quantities  (14) are simply the matrices formed from the
components  of  the  corresponding  curvature tensor $(cf. [6, 3, 4])$  and
that the set $\gamma (J^{p})$ consists of one or more neighborhoods in M.
Consequently,  now theorem 4 states that the investigated here special
bases exist iff the corresponding derivation is flat, i.e. if  its curvature
tensor is zero, a result already found in [6].

	In the general case, when  $2\le p<n ($if $n\ge 3)$, special bases,
even  anholonomic,  of the considered here type do not exist if (and only
if)  the  conditions  (17) are not satisfied. Besides, in this case  the
quantities (14) cannot  be considered as a "curvature" of $\gamma (J^{p})$
as  they  are something like "commutators" of covariant derivatives  of  a
type  $\nabla _{F}$,  where $F$ is a tangent to $\gamma (J^{p})$ vector field
(i.e.  $F\in T(\gamma (J^{p}))$  if  $\gamma (J^{p})$ is a submanifold of
$M)$, and which act on tangent to $M$ vector fields.

	Let us also note that  the bases in which the components of some
derivation vanish on a set $U$ are generally anholonomic, if any, and only in
the torsion free case (the derivation's torsion vanishes on $U)$ they may be
holonomic.

 	The  above results outline the general bounds of validity and  are
the exact mathematical expression of the equivalence principle,  which
states  [1] that the gravitational field strength, theoretically identified
 with  the components of a linear connection, can locally  be transformed  to
 zero by a suitable choice of the local  reference  frame  (basis),  i.e. it
 requires the existence of local bases in which the corresponding
 connection's components vanish.

 	 The above discussion, as well as the results from $[6, 7]$, show the
identical validity of the equivalence principle in zero and one dimensional
 cases, i.e. for $p=0$ and $p=1$. Besides these are the only  cases  when  it
 is fulfilled for arbitrary gravitational fields. In  fact,  for  $p\ge 2$
 (for  $n\ge 2$),  as we saw in Sect. IV, bases with the  above  property
 do not exist unless the conditions (17) are satisfied.  In  particular,
 for  $p=n\ge 2$ it is valid only for flat linear connections $(cf. [6])$.

 	Mathematically the equivalence principle is expressed through
 corollary~10 (or, in some more general situations, through theorem 4).
 Thus  we  see that in gravity theories based on linear connections this
 principle is identically satisfied at any fixed point or  along any  fixed
 path, but on submanifolds of dimension greater or  equal  two  it  is
 generally not valid. Therefore in this class of  gravity  theories  the
 equivalence principle is a theorem derived from their mathematical
 background. It may play a role as a principle if one tries to construct a
 gravity theory based on more general derivations, but, generally, it will
 reduce this theory to one based on linear connections.

\medskip
\medskip
 {\bf ACKNOWLEDGEMENT}
\nopagebreak

\medskip
 This  research  was partially supported by the Foundation for  Scientific
 Research of Bulgaria under contract Grant No. $F 103$.

\medskip
\medskip
 {\bf REFERENCES}

\medskip
 1. C.  W.  Misner,  Thorne K. S., Wheeler J. A., Gravitation (W. H.
 Freeman and Company, San Francisco, 1973). \par
 2. Schouten  J.  A.,  Ricci-Calculus:  An  Introduction  to  Tensor
 Analysis  and  its  Geometrical  Applications  (Springer Verlag,
 Berlin-G\"ottingen-Heidelberg, 1954), 2-nd ed. \par
 3. Kobayashi  S.,  K.  Nomizu, Foundations of Differential Geometry
 (Interscience Publishers, New York-London, 1963), vol. I. \par
 4. Lovelock D., H. Rund, Tensors, Differential Forms, and Variational
Principles (Wiley-Interscience Publications, John Wiley \& Sons, New
York-London-Sydney-Toronto, 1975).\par
 5. Rashevskii  P.  K.,  Riemannian  Geometry  and  Tensor  Analysis
 (Nauka, Moscow, 1967) (In Russian). \par
 6. Iliev  B.  Z., Special bases for derivations of tensor algebras.
 I. Cases in a neighborhood and at a point, Communication JINR, $E5-92-507$,
	Dubna, 1992.\par
7. Iliev B. Z., Special bases for derivations of  tensor  algebras. II. Case along paths, Communication JINR, $E5-92-508$, Dubna, $1992. 8$.    Raifeartaigh L., Fermi coordinates, Proceedings of  the  Royal Irish Academy, vol. {\bf 59}, Sec. {\bf A}, No. $2, 1958$. \par
 9. Hartman Ph., Ordinary Differential Equations (John Wiley \& Sons,
 New York-London-Sydney, 1964).

\end{document}